\documentclass{IEEEtran}

\usepackage{graphicx}
\usepackage{algorithm}
\usepackage{algorithmic}
\usepackage{booktabs}
\usepackage{hyperref}
\usepackage{enumitem}
\usepackage{wrapfig}
\usepackage{pifont}
\usepackage{amsmath, amsthm, amssymb, amsfonts, mathtools}

\usepackage[english]{babel}

\graphicspath{{Figures/}}

\def\R{\mathbf{R}}

\def\1{\mathbf{1}}
\def\conv{\mathrm{conv}}

\def\conv{\mathrm{conv}}

\newtheorem{theorem}{Theorem}
\newtheorem{lemma}{Lemma}
\newtheorem{corollary}{Corollary}

\newtheorem{observation}{Observation}

\title{Birkhoff's Decomposition Revisited: \\Sparse Scheduling for High-Speed Circuit Switches}
\author{
\IEEEauthorblockN{V\'ictor Valls\IEEEauthorrefmark{1}, George Iosifidis\IEEEauthorrefmark{2}, Leandros Tassiulas\IEEEauthorrefmark{1}} \\
\IEEEauthorblockA{\IEEEauthorrefmark{1}Yale University}
\IEEEauthorblockA{\IEEEauthorrefmark{2}TU Delft}
}
\begin{document}
\maketitle

\begin{abstract}
Data centers are increasingly using high-speed circuit switches to cope with the growing demand and reduce operational costs. One of the fundamental tasks of circuit switches is to compute a sparse collection of switching configurations to support a traffic demand matrix. Such a problem has been addressed  in the literature with variations of the approach proposed by Birkhoff in 1946 to decompose a doubly stochastic matrix exactly. However, the existing methods are heuristic and do not have theoretical guarantees on how well a collection of switching configurations (i.e., permutations) can approximate a traffic matrix (i.e., a scaled doubly stochastic matrix).

In this paper, we revisit Birkhoff's approach and make three contributions. First, we establish the first theoretical bound on the sparsity of Birkhoff's algorithm (i.e., the number of switching configurations necessary to approximate a traffic matrix). In particular, we show that by using a subset of the admissible permutation matrices, Birkhoff's algorithm obtains an $\epsilon$-approximate decomposition with at most $O( \log(1 / \epsilon))$ permutations. Second, we propose a new algorithm, \texttt{Birkhoff+}, which combines the wealth of Frank-Wolfe with Birkhoff's approach to obtain sparse decompositions in a fast manner. And third, we evaluate the performance of the proposed algorithm numerically and study how this affects the performance of a circuit switch. Our results show that \texttt{Birkhoff+}  is superior to previous algorithms in terms of throughput, running time, and number of switching configurations. 
\end{abstract}

\section{Introduction}

Data centers are increasingly adopting hybrid switching designs that combine traditional electronic packet switches with high-speed circuit switches \cite{FPR+10, CSS13,HSN16}.  In short, packet switches are flexible at making forwarding decisions at a packet level, but have limited capacity and are becoming increasingly expensive in terms of cost, heat, and power. In contrast, circuit switches provide significantly higher data rates at a lower cost but are less flexible at making forwarding decisions. The main drawback of circuit switches is that they have high reconfiguration times, which limit the amount of traffic they can carry \cite{WAK+10,FPR+10, LLF+14}. For instance, circuit switches have reconfiguration times in the order of milliseconds (e.g., $25$ ms for off-the-self circuit switches \cite{Polatis7000, Calient160}), whereas the reconfiguration times in electronic switches are in the scale of microseconds. As a result, hybrid switching architectures load balance and use circuit switches for high-intensity/bursty flows \cite{LML+15,BAV16}\footnote{Traffic in data centers is often bursty \cite{BAM10, RZB+15} and uses few input/output ports \cite{KSV+13}.}  and electronic switches for traffic that needs of a more fine-grained scheduling (e.g., delay-sensitive applications).

The problem of computing switching configurations for circuit switches is central to networking and has a direct impact on the performance of nowadays data centers. 
Mathematically, we can model a circuit switch as a crossbar,\footnote{See, for example, \cite[Section 4.1]{SY13}.} and cast the problem of finding a small collection of switching configurations as decomposing a doubly stochastic matrix\footnote{A matrix is doubly stochastic if its entries are non-negative and the sum of every row and column is equal to one. A permutation matrix is a binary doubly stochastic matrix.} as a \emph{sparse} convex combination of permutations matrices. In brief, for a given $n\times n$ doubly stochastic matrix $X^\star$ (i.e., a scaled traffic matrix) and an $\epsilon \ge 0$, the goal is to find a \emph{small} collection of permutation matrices $P_1,P_2,\dots, P_k$ (i.e., switching configurations) and weights $\theta_1,\theta_2,\dots,\theta_k >0$ (i.e., the fraction of time the switching configurations will be used) with $\sum_{i=1}^k \theta_i \le 1$ such that 
\begin{align}
\left \| X^\star - \sum_{i=1}^k \theta_i P_i \right \|_F \le \epsilon, \label{eq:ABD}
\end{align}
where $\| \cdot \|_F$ is the Frobenius norm (see definition in Section \ref{sec:notation}). The smaller $\epsilon$ is, the higher the throughput. However, practical systems cannot use as many switching configurations as desired as each inflicts a reconfiguration time $\delta$ that affects the fraction of time the switch can carry traffic.\footnote{Technically, a traffic matrix $X^\star$ is valid for a time window period $W$, and the decomposition must satisfy $\sum_{i=1}^k (\theta_i + \delta) \le W$. That is, the time spent transmitting $(\sum_{i=1}^k \theta_i)$ and reconfiguring ($\delta k$) cannot exceed the time window duration ($W$). } Or put differently, there is a constraint on the number of configurations a switch can use to approximate a traffic matrix.

Previous work has addressed the problem above with variations (e.g., \cite{CCH00,LML+15,DU16,LV18}) of the approach proposed by Birkhoff in 1946 \cite{Bir46} to decompose a doubly stochastic matrix \emph{exactly} (i.e., $\epsilon = 0$). However, little is known about the behavior or convergence properties of Birkhoff's algorithm, and so fundamental questions remain still unanswered. In particular, how does the decomposition approximation $\epsilon$ in Eq.\ \eqref{eq:ABD}  depend on the number of switching configurations?  How much does an additional switching configuration contribute to increasing a circuit switch throughput? How is Birkhoff's algorithm related to other numerical methods in other fields, such as optimization and machine learning? Answering these questions is crucial to better understand the structure of the problem and to design new algorithms that improve the performance of circuit switches. To this end, the main contributions of the paper are the following:

\begin{table}
\caption{Sparsity and permutation selection complexity of \texttt{Birkhoff+} (this paper) and previous algorithms. LP and QP stand for linear and quadratic program respectively. }

\begin{tabular}{lllll}
\toprule
Algorithm & Sparsity ($k$) & Perm. selec. complexity  \\
\midrule
\texttt{Birkhoff} \cite{Bir46}&  ---  & One LP   \\
\texttt{Solstice}\cite{LML+15}  & ---  & Multiple LPs  \\
\texttt{Eclipse}$^*$ \cite{BAV16}  &  Approx. ratio  & Multiple LPs  \\ 
\texttt{FW} \cite{FW56} &  $O(1/\epsilon^2)$ & One LP   \\
\texttt{FCFW} \cite{LJ15}&  $O(\log(1 /\epsilon^2))$   & One LP + QP($k$) \\
\texttt{Birkhoff+} (this paper) &  $O(\log( 1 /\epsilon))$  & One LP \\
\bottomrule
\end{tabular}
\label{table:algorithms}
\end{table}

\textbf{(i) Revisiting Birkhoff's approach.} We revisit Birkhoff's algorithm and establish the first theoretical bound on its sparsity (i.e., the number of permutations necessary to approximate a doubly stochastic matrix).
In particular, we show that by selecting permutations from a subset of admissible permutations, Birkhoff's algorithm has sparsity $O(\log(1 / \epsilon))$ (Theorem \ref{th:birkhoff_linear_convergence}). That is, the  number of permutations required to obtain an $\epsilon$-approximate decomposition increases logarithmically with the decomposition error. Our results also show that previous Birkhoff-based algorithms that select permutations using a Max-Min criterion (e.g., \cite{LML+15}) have logarithmic sparsity (Corollary \ref{th:coro_greedy}), and that Birkhoff's algorithm is strongly connected to block-coordinate descent and Frank-Wolfe methods in convex optimization (Section \ref{sec:block-coordinate} and Section \ref{sec:frank-wolfe}).

\textbf{(ii) New algorithm (\texttt{Birkhoff+}).} We propose a new algorithm that combines Birkhoff's approach and Frank-Wolfe. Specifically, permutation matrices are selected using a Frank-Wolfe-type update with a barrier function, while the weights as in Birkhoff's approach. The proposed algorithm has theoretical guarantees (Corollary \ref{th:bir_coro}) and is non-trivial as a direct combination of Birkhoff's approach with Frank-Wolfe may not converge (Theorem \ref{th:birkhoff_algorithm_fwform_notcvg}). 
Furthermore, \texttt{Birkhoff+} is faster than previous algorithms as it computes a new permutation/configuration by solving a \emph{single} linear program (LP). Table \ref{table:algorithms} contains a summary of the main differences between \texttt{Birkhoff+} and the state-of-the-art algorithms discussed in Section \ref{sec:related_work}.

  \textbf{(iii) Numerical evaluation.} We evaluate \texttt{Birkhoff+}'s performance in a circuit switch application  and compare it against existing algorithms for a range of matrices (dense, sparse, skewed) that capture the characteristics of traffic in data centers. Our results show that \texttt{Birkhoff+}  is superior to previous algorithms in terms of throughput, running time, and number of switching configurations. For instance, when $\delta / W = 10^{-2}$ (the reconfiguration time over the time available for transmission), \texttt{Birkhoff+} has $7\%$ more throughput than the best state-of-the-art algorithm. If we consider, in addition, the time to compute the switching configurations as an overhead, the throughput gain increases  to $34\%$ (switch with $n =100$ ports).

The outline of the paper is as follows. Section \ref{sec:related_work} presents related work and Section \ref{sec:preliminaries} the preliminaries, which include the notation and how to find a permutation matrix by solving a linear program. In Section \ref{sec:birkhoff_revisited}, we revisit Birkhoff's approach in a \emph{general} form, establish its sparsity rate, and show how this is connected to block-coordinate descent methods in convex optimization. The latter also clarifies that selecting a permutation matrix can be seen as choosing a (gradient) descent direction. 
Section \ref{sec:frank-wolfe} shows how to use Frank-Wolfe algorithms to decompose a doubly stochastic matrix, and how Frank-Wolfe chooses a permutation matrix that provides ``steepest descent.'' In Section \ref{sec:birkhoff+}, we present the new algorithm (\texttt{Birkhoff+}) and in Section \ref{sec:numerical_evaluation} evaluate its performance against the state-of-the-art algorithms. Finally, Section \ref{sec:conclusions} concludes. All the proofs are in the Appendix.

\section{History and related work } 
\label{sec:related_work}

\subsection{Birkhoff's approach} This is the method employed by Birkhoff in 1946 to decompose a doubly stochastic matrix \emph{exactly} \cite[first theorem]{Bir46}.\footnote{The result is also known as Birkhoff-von Neumman (BvN) as it was discovered independently by von Neumman \cite{Neu53}. We use Birkhoff instead of BvN as the algorithm used in the literature is based on the method of proof used by Birkhoff in \cite{Bir46}.} In brief, the method consists of finding permutations matrices sequentially (e.g., with the Hungarian algorithm) and terminates when it obtains an exact decomposition, which happens with at most $k = (n-1)^2 + 1$ iterations/permutations \cite{Mar60,Bru82} by Carath\'eodory's theorem.
The method, however, does not guarantee that the decomposition (i.e., $\sum_{i=1}^k \theta_i P_i$) is close to the doubly stochastic matrix it aims to approximate (i.e., $X^\star$). In fact, the approximation is typically very poor until the algorithm converges exactly in the last iteration (see Figure \ref{fig:sparsity}a in Section \ref{sec:numerical_evaluation} for an example).

\subsection{Related mathematical problems} The problem of finding the Birkhoff decomposition with the minimum number of permutation matrices ($\min k$ s.t. $X^\star = \sum_{i=1}^k \theta_i P_i$) was addressed in \cite{DU16} and shown to be NP-hard. In \cite{KLS17}, the authors also show that the problem is not tractable when the minimal decomposition can be expressed with $k \ge 4$ permutations. 
The work in \cite{LH03} formulates a similar problem. For a demand matrix $D =  X^\star - S$ with $S \in [0,1]^{n \times n}$,\footnote{The entries of the demand matrix $D$ are non-negative. Matrix $S$ adds a non-negative virtual load to demand matrix so that $D + S$ is doubly stochastic.  } 
the goal is to find a collection of weights $\{\theta_i \}_{i=1}^k$  and permutation matrices $\{ P_i \}_{i=1}^k$ that minimizes $\sum_{i=1}^k (\theta_i + \delta )$ subject to $\sum_{i=1}^k \theta_i P_i \ge D$ entry-wise. The problem is shown to be NP-complete. The problem addressed in this paper is different in spirit from the mathematical problems in \cite{DU16,LH03} because we do not aim to find a (small) collection of objects (i.e., $k$) subject to decomposition constraints (i.e., $X^\star = \sum_{i=1}^k \theta_i P_i$ or $\sum_{i=1}^k \theta_i P_i \ge D$). Instead, our goal is to design an algorithm that minimizes $\| X^\star - X_k \|_F$ where $X_k = \sum_{i=1}^k \theta_i P_i$. The convergence rate of the numerical method correponds the number of permutations required to obtain an $\epsilon$-approximate decomposition.  

\subsection{Algorithms} The paper in \cite{LML+15} proposes \texttt{Solstice}, a  Birkhoff-based heuristic for finding a Birkhoff decomposition with few permutations/configurations. \texttt{Solstice} picks permutation matrices using a Max-Min type criterion, and the weights or configurations durations are selected as large as possible provided $X^\star - \sum_{i=1}^k \theta_i P_i$ is non-negative entry-wise.  
The work in \cite{BAV16} proposes \texttt{Eclipse}, a sub-modular-type algorithm for solving the problem of the type introduced in \cite{LH03}. Permutation matrices and weights are selected jointly to maximize an effective utilization criterion, which takes into account the reconfiguration penalty $\delta$. Also, \cite{BAV16} shows that the final decomposition satisfies the optimal approximation ratio in sub-modular optimization with cover constraints. 
Both algorithms  \cite{LML+15, BAV16} select permutation matrices by solving multiple linear programs (LPs) with a simplex type method \cite{Dan63}. Finally, we note the recent works in \cite{LV18} and \cite{SSY19}. The first extends \texttt{Eclipse} to use a special type of weights/time coefficients that do not constraint the decomposition to be a scaled doubly stochastic matrix. The second addresses the \emph{online} version of the problem in \cite{BAV16}---in the machine learning sense \cite{Bub11}---where the traffic matrix is learned a posteriori.

To conclude, we note the Frank-Wolfe algorithms \cite{FW56,Jag13} used extensively in machine learning. The Frank-Wolfe setup is the following. Given a collection of discrete objects $\mathcal D$ and a convex set $\mathcal X \subseteq \mathrm{conv}{(\mathcal D)}$,  the goal is to minimize a convex function  by making convex combinations of the discrete objects. The problem addressed in this paper can be seen as a special case for Frank-Wolfe. The permutation matrices correspond to the discrete objects, the Birkhoff polytope is the convex set, and the objective function a metric that captures the distance between the approximate decomposition and $X^\star$ (e.g., Frobenius norm or Euclidean distance). Also, and unlike Birkhoff-based approaches, Frank-Wolfe algorithms provide sparsity guarantees and ensure that the approximate decomposition is always a doubly stochastic matrix.


\section{Preliminaries}
\label{sec:preliminaries}

\subsection{Notation} 
\label{sec:notation}
We use $\R_+$ and $\R^d$ to denote the set of nonnegative real numbers and $d$-dimensional real vectors.  Vectors and matrices are written in lower and upper case respectively, and all vectors are in column form. The transpose of a vector $x \in \R^d$ is indicated with $x^T $, and we use $\1$ to indicate the all ones vector---the dimension of the vector will be clear from the context. We use parenthesis to indicate an element in a vector, i.e., $x{(j)}$ is the $j$'th element of vector $x$. Similarly, the element in the $i$'th row and $j$'th column of a matrix $X$ is indicated with $X(i,j)$.
For two vectors $x, y  \in \R^d$, we write $x \succ y$ when $x{(j)} > y{(j)}$ for all $j \in \{ 1,\dots,d\}$, and $x \succeq y$ when $x(j) \ge  y(j)$.  Finally, we recall the Frobenius norm of a matrix $X$ is defined as $\| X \|_F = \sqrt{\sum_{i,j} | X(i,j)|^2} = \sqrt {\mathrm{Tr} ( X X^*)}$ and that $[n]$ is the  short-hand notation for $\{1,\dots,n\}$.

\subsection{Finding extreme points by solving linear programs}
\label{sec:findingextremepoints}
We will present algorithms that find extreme points (i.e., permutation matrices) by solving linear programs (LPs) over a convex set (i.e., the Birkhoff polytope or set of doubly stochastic matrices). We recall the following result from linear programming.

\begin{lemma}
Let $\mathcal X$ be a bounded polytope from $\R^d$, and $\mathcal E$ denote its extreme points. For any vector $c \in \R^d$, we have that $\{ \arg \min_{x \in  \mathcal X} \ c^T x \} \cap  \mathcal E \ne \emptyset$.
\label{th:LPep}
\end{lemma}
That is, an extreme point in $\mathcal E$ is always a solution to $\min_{x \in \mathcal X} c^T x$. 
 In our case, $\mathcal X$ is the Birkhoff polytope and $\mathcal P$ the set of permutation matrices. 
Throughout the paper, we will cast linear programs as 
\begin{align}
\texttt{LP}(c, \mathcal X):  
\begin{array}{llll}
& \text{minimize} & c^T x \\
&  \text{subject to} & x \in \mathcal X
\end{array},
\end{align}
and we will assume that the solution returned is always an extreme point---which is the case if we solve the LP with a simplex-type method \cite{Dan63}.

\section{Revisiting Birkhoff's Algorithm}
\label{sec:birkhoff_revisited}

This section revisits Birkhoff's algorithm. The main technical contribution is Theorem \ref{th:birkhoff_linear_convergence}, which establishes that the number of permutation matrices in Birkhoff's approach increases logarithmically with the decomposition error. 

\subsection{Approximate Birkhoff decomposition problem }

The mathematical problem we want to solve is the following. For a given $n \times n$ doubly stochastic matrix $X^\star$ and an $\epsilon \ge 0$, our goal is to find a \emph{small} collection of permutation matrices $P_1,P_2,\dots, P_k$ and weights $\theta_1,\theta_2,\dots,\theta_k >0$ with $\sum_{i=1}^k \theta_i \le 1$ such that $\| X^\star - \sum_{i=1}^k \theta_i P_i \|_F \le \epsilon$. Recall a matrix $X \in [0,1]^{n \times n}$ is doubly stochastic if every row and column sums to one. That is, $X \1 = \1$ and $\1^T X = \1^T$. Also, a doubly stochastic matrix is a permutation if its entries are binary.

\subsection{Algorithm description}

The original Birkhoff algorithm is described in Algorithm \ref{al:birkhoff_algorithm}, and consists of two steps. First, the algorithm calls the subroutine \texttt{PERM}, which returns a permutation $P_k$ and a weight $\theta_k$. The second step is to add $\theta_k P_k$ to the previous approximate decomposition, i.e., $X_k = X_{k-1} + \theta_k P_k$. 
The permutation $P_k$ and weight $\theta_k$ must satisfy the following three conditions:
\begin{align}
& \textstyle X_{k-1}(a,b) + \theta_k P_k (a,b) \le X^\star(a,b) && \forall a,b \in[n] \label{eq:k1} \\
& \textstyle \theta_k > 0 && \forall k \ge 1 \label{eq:k2} \\
& \textstyle \sum_{i=1}^k \theta_i \le 1 && \forall k\ge 1  \label{eq:k3}
\end{align}
In words, $X_k(a,b) \le X^\star (a,b)$ for all $a,b \in \{1,\dots,n\}$, the weights are strictly positive, and the sum of the weights is less than or equal to one.  The algorithm terminates when the approximate decomposition $X_k$ is $\epsilon$ close to $X^\star$, or when the maximum number of admissible permutations ($k_\text{max}$) is reached. 

\begin{algorithm}[t]
\caption{\texttt{General Birkhoff}}
\begin{algorithmic}[0]
\STATE \textbf{Input:} Doubly stochastic matrix $X^\star$, $\epsilon \ge 0$, and $k_\text{max} \ge 1$
\STATE \textbf{Set:} $k = 1$ and $X_0 =\{0 \}^{n \times n}$
\WHILE{$\| X_{k-1} - X^\star \|_F  > \epsilon $ and $k \le k_\text{max}$}
\STATE $P_k, \theta_k \leftarrow \texttt{PERM}(X_{k-1},X^\star)$
\STATE $X_{k} \leftarrow X_{k-1} + \theta_k P_k$
\STATE $k \leftarrow k + 1$ 
\ENDWHILE
\RETURN $(P_1,\dots, P_{k-1})$, $(\theta_1,\dots,\theta_{k-1})$
\end{algorithmic}
\label{al:birkhoff_algorithm}
\end{algorithm}

\subsection{Convergence}
\label{sec:convergence}

We proceed to establish the convergence of Birkhoff's algorithm. We start by presenting the following lemma, which establishes a lower and upper bound on $\| X_k - X^\star \|_F$.

\begin{lemma} 
\label{th:birkhoff_lemma}
Consider the setup in Algorithm \ref{al:birkhoff_algorithm} and suppose the subroutine 
\texttt{PERM} returns a weight $\theta_k$ and a permutation matrix $P_k$ that satisfy the conditions in Eqs.\ \eqref{eq:k1}- \!\eqref{eq:k3} for all $k \ge 1$. Then, the following two bounds hold:
\begin{align}
 & \| X_{k} - X^\star \|_F \ge \left(1 - \sum_{i=1}^k \theta_i\right)  \label{eq:birkhoff_lower_bound} \\
 & \| X_{k} - X^\star \|_F  \le  \sqrt{ n  \prod_{i=1}^k  \left(1  - \frac{n \theta^2_i}{\| X_{i-1} - X^\star \|_F^2} \right) }
\label{eq:birkhoff_upper_bound}
\end{align}
where $\theta_i \le \frac{1}{\sqrt n}{\| X_{i-1} - X^\star \|_F}$.

\end{lemma}

The bounds in Lemma \ref{th:birkhoff_lemma} are very general as they hold for any collection of permutation matrices and weights that satisfy the conditions in  Eqs.\ \eqref{eq:k1}-\eqref{eq:k3}. 
The lower bound in Eq.\ \eqref{eq:birkhoff_lower_bound} tells us that the approximate decomposition error is at least $(1 - \sum_{i=1}^k \theta_i)$, and so we will have an exact decomposition (i.e., serve $100\%$ of the traffic demand) only if $\sum_{i=1}^k \theta_i = 1$.
The upper bound in Eq.\ \eqref{eq:birkhoff_upper_bound} shows how the decomposition error depends on the weights $\theta_i$ and the previous approximations  $\| X_{i-1} - X^\star \|_F^2$, $i=1,\dots,k$. In particular, on the ratio ${n \theta^2_i}/{\| X_{i-1} - X^\star \|_F^2}$, which captures how large $\theta_i$ is with respect to the previous approximation. Note that the values that $\theta_i$ can take depend on $\| X_{i-1} - X^\star \|_F^2$ as we must always satisfy the conditions in Eqs.\ \eqref{eq:k1}--\eqref{eq:k3}. Finally, we note that finding a joint collection of weights and permutation matrices that minimize the RHS of Eq.\ \eqref{eq:birkhoff_upper_bound} for a fixed $k$ is as difficult as minimizing $\| X_k - X^\star\|_F^2$ directly, since the RHS of Eq.\ \eqref{eq:birkhoff_upper_bound} depends on $\| X_{i-1} - X^\star \|_F$, $i=1,\dots,k$. Because of the latter, we study how to minimize the RHS of Eq.\ \eqref{eq:birkhoff_upper_bound} in an iterative manner: for a given collection of permutation matrices $P_i$ and weights $\theta_i$ with $i=1,\dots,k-1$, our goal is to find a permutation matrix $P_k$ and weight $\theta_k$ that decrease the RHS of Eq.\ \eqref{eq:birkhoff_upper_bound}.

In the following, we study the algorithm's progress in terms of error for every additional permutation matrix in the decomposition.  Addressing this question is important to obtain a bound on the number of permutations required to obtain an $\epsilon$-approximate decomposition as well as to know how to select good permutation matrices.  To start, let $\mu_i \coloneqq n \theta^2_i/ \| X_{i-1} - X^\star \|_F^2$ and rewrite Eq.\ \eqref{eq:birkhoff_upper_bound} as
\begin{align}
\| X_{k} - X^\star \|_F  \le  \sqrt{ n  \prod_{i=1}^k  (1-\mu_i) } \label{eq:bound_mu}
\end{align}
Note that $(1-\mu_i) \in [0, 1)$ for all $i=1,2,\dots$ since $\theta_i \le \frac{1}{\sqrt n}{\| X_{i-1} - X^\star \|_F}$. Hence, we have that  $X_{k} \to  X^\star$ as $k \to \infty$ and so the algorithm converges. Now, suppose there exists a constant $\mu_{\min} > 0$ such that $\mu_{\min} \le  \mu_i$ for all $i \ge 1$. Then, the bound in Eq.\ \eqref{eq:bound_mu} simplifies to
\begin{align}
\| X_{k} - X^\star \|_F  \le  \sqrt n   (1-\mu_{\min})^{k/2}. \label{eq:bir_lin_cvg}
\end{align}
The last equation tells us that the approximation error decreases exponentially with the number of permutations. For example, if $ \mu_{\min} = 1/2$, we have that $\| X_{k} - X^\star \|_F  \le  \sqrt n  (1/2)^{k/2}$, which means that every additional permutation in the decomposition decreases the approximation error by at least half. The ratio $\kappa \coloneqq 1/\mu_{\min} \ge 1$ can be regarded as the condition number in optimization with a strongly convex objective \cite{BV04}[Section 9.1.2 and 9.3.1]. 

The following lemma establishes an upper bound on the number of permutations required to obtain an $\epsilon$-approximate decomposition provided that a constant $\mu_{\min} > 0$ exists. 

\begin{lemma}
\label{th:lemma_decomposition}
Suppose ${n \theta^2_i}/{\| X_{i-1} - X^\star \|_F^2} \ge \mu_{\min}$ for all $i=1,\dots,k$ for some constant $\mu_{\min} > 0 $. Then, Algorithm \ref{al:birkhoff_algorithm} obtains an $\epsilon$-approximate decomposition with at most
\[
k \le 2 \log^{-1} \left(\frac{1}{1-\mu_{\min}} \right) \log \left( \frac{\sqrt n}{\epsilon}  \right)
\]
permutation matrices.
\end{lemma}

Lemma \ref{th:lemma_decomposition} says that if a constant $\mu_{\min}$ exists, then the number of permutation required to obtain an $\epsilon$-approximate decomposition has a logarithmic dependence with $\epsilon$. Hence, it remains to show whether such constant exists. Or equivalently, we need to show that we can select a $\theta_i$ such that  ${n \theta^2_i}/{\| X_{i-1} - X^\star \|_F^2}$ is uniformly lower bounded by a strictly positive constant. We show that in the following theorem, which is one of the main contributions of the paper.

\begin{algorithm}[t]
\caption{Subroutine \texttt{PERM}}
\begin{algorithmic}[1]
\STATE \textbf{Input:} $X^\star$ and $X_{k-1} = \sum_{i=1}^{k-1} \theta_i P_i$
\STATE  $\alpha \leftarrow (1-\sum_{i=1}^{k-1} \theta_i ) / n^2 $
\STATE  $P_k \leftarrow \hat P \in \mathcal I_k (\alpha)$
\STATE  $\theta_k \leftarrow \texttt{BIRKHOFF\_STEP}(X^\star, X_{k-1}, P_k)$ \ \  (Algorithm \ref{al:birkhoff_step})
\end{algorithmic}
\label{al:birkhoff_weights}
\end{algorithm}

\begin{algorithm}[t]
\caption{\texttt{BIRKHOFF\_STEP}}
\begin{algorithmic}[1]
\STATE \textbf{Input:} $X^\star$, $X_{k-1}$, and $P_k$
\RETURN $\min_{a,b} \left\{ (X^\star(a,b) - X_{k-1}(a,b) - 1 ) P_k(a,b)  + 1 \right\}$ 
\end{algorithmic}
\label{al:birkhoff_step}
\end{algorithm}

\begin{theorem}
\label{th:birkhoff_linear_convergence}
Let $\mathcal P$ be the set of $n\times n$ permutation matrices and define 
\begin{align}
\! \mathcal I_k (\alpha) = \left\{ P \in \mathcal P \mid   
X_{k-1}(a,b)  + \alpha P(a,b)  \le X^\star(a,b) \right\}
\label{eq:basic_penalty}
\end{align}
with $\alpha = \sqrt{\frac{\mu_{\min}}{ n} } (1-\sum_{i=1}^k \theta_i)$ where $\mu_{\min} = 1 / n^3$. 
Then, $\mathcal I_k (\alpha) \ne \emptyset$, and Algorithm \ref{al:birkhoff_algorithm} with the subroutine \texttt{PERM} defined in Algorithm \ref{al:birkhoff_weights} obtains an $\epsilon$-approximate decomposition with at most
\begin{align}
k \le 2 \log^{-1} \left(1 - \underset{i \in [k]}{\min} \frac{n \theta^2_i}{\| X_{i-1} - X^\star \|_F^2} \right)^{-1}  \log \left( \frac{\sqrt n}{\epsilon}  \right) \label{eq:main_bound}
\end{align}
permutation matrices.
\end{theorem}

Theorem \ref{th:birkhoff_linear_convergence} establishes that by selecting permutation matrices from set $\mathcal I_k(\alpha) \subseteq \mathcal P$, and weights as indicated in Algorithm \ref{al:birkhoff_weights}, then the number of permutation matrices required to obtain an $\epsilon$-decomposition increases logarithmically with $\epsilon$. 
Set $\mathcal I_k(\alpha)$ is necessary to enforce that conditions in Eqs.\ \eqref{eq:k1}--\eqref{eq:k3} are satisfied, but also to push Birkhoff's algorithm to make sufficient progress in every iteration. Observe that the threshold $\alpha$ is bounded away from zero and that this depends on the constant $\mu_{\min} = 1 / n^3$. 
Finally, we have written $ \min_{i \in \{1,\dots,k\}} {n \theta^2_i}/{\| X_{i-1} - X^\star \|_F^2}$ instead of $\mu_{\min}$ in Eq.\ \eqref{eq:main_bound} (c.f.\ Lemma \ref{th:lemma_decomposition}) to emphasize two points. The first one is that $\mu_{\min}$ is over-conservatively small, and that we can in general obtain a much sharper upper bound. In the numerical evaluation (Section \ref{sec:zero_exp}), we show the condition numbers ($\kappa = 1 / \mu_{\min}$) of different algorithms.  
The second point is that $n \theta^2_i / \| X_{i-1} - X^\star \|_F^2$ is a quantity that we can measure and so use a as a criterion for selecting a ``good enough'' permutation matrix. Importantly, recall that the \texttt{PERM} subroutine does not specify which specific permutation to select from $\mathcal I_k(\alpha)$, which is in marked contrast to pervious approaches (e.g., \cite{LML+15,BAV16,LV18}), which use a predefined criterion for selecting permutation matrices and weights.

\subsection{Discussion}

\subsubsection{Max-Min Birkhoff algorithms}

The most popular variant of Birkhoff's algorithm (e.g.,  \cite{LML+15,DU16}) aims to find a permutation matrix with the largest associated weight. Such approach corresponds to solving the following optimization problem:
\begin{align}
\begin{array}{lllll}
\underset{\theta > 0, P \in \mathcal P}{\text{maximize}} & \theta \\
\text{subject to} &  X_{k-1}(a,b) + \theta P(a,b) \le X^\star(a,b)  \\& \forall a,b \in [n]
\end{array}
\label{eq:weight_opt_problem}
\end{align}
The strategy is also known as Max-Min because it is equivalent to finding a permutation matrix $P$ with the largest smallest element $X^\star(a,b) - X_k(a,b)  $ provided $P(a,b) = 1$. Hence, the set of solutions to the optimization problem above is given by
\[
\mathcal S_k \coloneqq \arg \max_{P \in \mathcal P} \left\{ \min_{\substack{a,b \in [n] \\ P(a,b) = 1} } X^\star(a,b) - X_k(a,b)  \right\}
\]
Note that $\mathcal S_k \subseteq \mathcal I_k(\alpha)$ since $\mathcal I_k(\alpha)$ includes all the solutions with $\theta \ge \alpha > 0$. Further, we have that $\mathcal I_k(\alpha) \ne \emptyset$ by Theorem \ref{th:birkhoff_linear_convergence}. We have arrived at the following corollary to Theorem \ref{th:birkhoff_linear_convergence}. 
\begin{corollary}[Theorem \ref{th:birkhoff_linear_convergence}]
\label{th:coro_greedy}
The Birkhoff-type algorithms that select permutation matrices using a Max-Min criterion (e.g., \cite{LML+15}) have sparsity $O( \log ( 1/{\epsilon}))$.
\end{corollary}

To conclude, we would like to emphasize that finding a permutation matrix in set $\mathcal S_k $ is non-trivial. The typical approach is to fix a weight $\theta$, and then try to find a permutation matrix that satisfies the constraints in Eq.\ \eqref{eq:weight_opt_problem}. The process is repeated for different weights, which are selected with different strategies; for example, \cite{LML+15} uses a halving threshold rule. The main issue with this method is that it is slow, and so non-convenient for applications that need to carry out decomposition fast. For example, when we are given a traffic matrix associated with a time window. The time spent computing the switching configurations is time that the switch cannot use for serving traffic. 

\subsubsection{Birkhoff's algorithm as a block-coordinate descent}
\label{sec:block-coordinate}
The Birkhoff algorithm can be thought in convex optimization terms. In particular, as solving the following convex optimization problem 
\begin{align}
\begin{tabular}{lll}
$\underset{X \in \R^{n\times n}}{\text{minimize}}$ & $\| X - X^\star \|_F^2$ \\
subject to & $X(a,b) \le X^\star (a,b)$ & $\forall a,b [n]$\\
 & $X(a,b) \ge 0$ & $\forall a,b \in [n]$ \\
\end{tabular}
\end{align}
using a block-coordinate descent method with $X_0 = \{0\}^{n \times n}$ (see \cite{Wri15} \cite{Nes12}, \cite{BT89}[Section 7.5.3]). Note that the objective is convex and the constraints linear. The block-coordinate method consists of the update\footnote{The method is usually expressed in vector form. In our case, we can create a vector by stacking the matrix columns.}
\[
X_k = X_{k-1} + \theta_k M_k
\]
where $\theta_k > 0$ is a step size and $M_k \in \{-1,0,1\}^{n\times n}$ a matrix that indicates the direction in which to update each of the coordinates. Birkhoff's approach can be regarded as a special case where the $M_k$ matrices are permutations, and so have constrains on the group of coordinates can be \emph{jointly} updated. Also, there are no negative coordinates since by selecting $X_0 = \{ 0 \}^{n\times n}$ as starting point the algorithm only needs to ``move forward.'' To conclude, we note that our sparsity result is connected to the linear convergence rate obtained by convex optimization algorithms that exploit the strong convexity of the objective function.

\section{Frank-Wolfe for the Approximate Birkhoff Decomposition}
\label{sec:frank-wolfe}

In this section, we show how the Frank-Wolfe (\texttt{FW}) algorithm and its fully corrective variant (\texttt{FCFW}) can be used to decompose a doubly stochastic matrix. The main contributions are to give explicit sparsity bounds for the \texttt{FW} and  \texttt{FCFW} algorithms (Theorem \ref{th:FW} and \ref{th:FCFW}) and to discuss the properties of how Frank-Wolfe selects permutation matrices (Observations \ref{ob:1} and \ref{ob:2}). The latter will be key to choose permutation matrices in the Birkhoff-type algorithm we will present in Section \ref{sec:birkhoff+}.

\begin{algorithm}[tb]
\caption{\texttt{Birkhoff (vector form)} }
\begin{algorithmic}[1]
\STATE \textbf{Input:} Birkhoff polytope $\mathcal B$, $x^\star \in \mathcal B$, $\epsilon \ge 0$, $k_\text{max} \ge 1$
\STATE \textbf{Set:} $k = 1$ and $x_0 = 0$
	\WHILE {$\| x^\star - x_{k-1} \|_2 >  \epsilon$ and $k \le k_\text{max}$}  
	\STATE  $\circ$ $ p_k \leftarrow$ \texttt{LP}  $(- \lceil x^\star - x_{k-1} \rceil, \mathcal B)$\\ 
	\STATE $\text{\footnotesize \ding{73}}$  $\theta_k \leftarrow \texttt{BIRKHOFF\textunderscore STEP}  (x^\star, x_{k-1}, p_k) $ \\ 
	\STATE $\diamond$ $x_{k} \leftarrow  x_{k-1} +  \theta_k p_k $ \\
	\STATE$ \phantom{a}$ $k \leftarrow k + 1$
	\ENDWHILE
 \RETURN $(p_1,\dots, p_{k-1})$, $(\theta_1,\dots,\theta_{k-1})$
\end{algorithmic}
\label{al:birkhoff_algorithm_fwform}
\end{algorithm}

\subsection{Birkhoff polytope and algorithm in vector form}
In the rest of the paper, it will be more convenient to write  $n \times n$ doubly stochastic matrices as  $n^2$-dimensional vectors\footnote{Instead of having a matrix $Z \in \R^{n\times n}_+$ such that $Z\1 = Z^T \1 = \1$, we work with a vector $x \coloneqq (z_1,\dots,z_n)$ where $z_i$ is the $i$'th column of $Z$.} in the set 
\[
\mathcal B \coloneqq \{ x \in \R^d \mid x \succeq 0, \ Ax= b\},
\]
where $d = n^2$, $A \in \{0,1\}^{2n \times d}$, and $b \coloneqq \{1\}^{2n}$. Matrix $A$ contains the $2n$ equality constraints that characterize the Birkhoff polytope (i.e., the sum of the columns and rows of a doubly stochastic matrix must be equal to $1$). 
The specific structure of $A$ can be derived easily and is given in the Appendix. 
As before, we use set $\mathcal P \subset \{0,1 \}^d$ to denote the set of permutation matrices or extreme points, but now these are in column form. The terms extreme point and permutation matrix will be used interchangeably in the rest of the paper. Finally, Algorithm \ref{al:birkhoff_algorithm_fwform} contains the procedure of the classic Birkhoff  algorithm \cite{Bir46} in vector form,\footnote{The algorithm corresponds to the method of proof employed by Birkhoff to show that a doubly stochastic matrix is an arithmetic measure of permutation matrices. See \cite{Bir46}, theorem on page 1.} which is a special case of the more general Algorithm \ref{al:birkhoff_algorithm}. Permutation matrices are selected by solving the liner program \texttt{LP}  $(- \lceil x^\star - x_{k-1} \rceil, \mathcal B)$ (see Section \ref{sec:findingextremepoints}) and the step sizes as large as possible provided $x_k \preceq x^\star$ for all $k\ge 1$.  The \texttt{LP}  $(- \lceil x^\star - x_{k-1} \rceil, \mathcal B)$ returns any admissible permutation matrix (see Section \ref{sec:preliminaries}) and $\lceil \cdot \rceil$ denotes the entry-wise ceiling of a vector.

\subsection{Frank-Wolfe overview}

In short, the Frank-Wolfe algorithm is a numerical method for minimizing a convex function $f$ over a convex set contained in the convex hull of a set of discrete points or atoms \cite{Jag13}. In our case, the convex set is the Birkhoff polytope ($\mathcal B$) and the atoms  the set of permutation matrices ($\mathcal P$). In each iteration, the algorithm selects an extreme point with update
\begin{align}
p_k \in \arg \min_{u \in \mathcal P } \  \nabla f(x_{k-1})^T u \label{eq:FWupdate}
\end{align}
and choses a step size $\theta > 0$ such that   $f(x_{k-1} + \theta (p_k-x_{k-1})) < f(x_{k-1})$. The essence of the algorithm is that when $f$ is smooth on $\mathcal B$,\footnote{There exists a constant $L$ such that $f(y) \le f(x) + \nabla f(x)^T (y-x) + \frac{L}{2} \| y-x \|_2^2$ for all $x,y \in \mathcal B$.} there always exists an extreme point that is a direction in which it is possible to improve the objective function. 
The step size can be selected in a variety of ways (e.g. constant, line search, etc.) and differently from the previous section, Frank-Wolfe does \emph{not} require that $x_{k-1} + \theta_k p_k \preceq x^\star$ where $x^\star$ is the doubly stochastic matrix we want to decompose. Also, Frank-Wolfe ensures, by construction, that $x_k$ is a convex combination of the permutation matrices throughout the iterations. As objective function, we use $f(x) = (1/2)\| x - x^\star\|_2^2$ to streamline exposition but also because it allows us to make the following observations:

\begin{figure}
  \begin{center}
 \includegraphics[width=0.40\columnwidth]{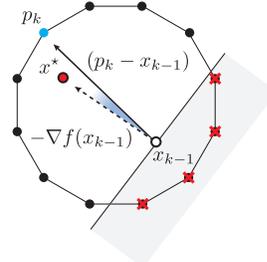}
  \end{center}
  \vspace{-20pt}
  \caption{Schematic illustration of the steepest descent permutation discussed in Observation \ref{ob:2}. The black dots with a red cross are the extreme points that are non-descent directions. Frank-Wolfe with $f(x)  = (1/2) \| x^\star - x\|^2_2$  chooses the extreme point $p_k$ (i.e., the permutation) that minimizes the angle between $(p_k - x_{k-1})$ and $-\nabla f(x_k) = (x^\star - x_{k-1})$. }
  \label{fig:steepest_descent}
\end{figure}

\begin{observation}[Weighted search direction] For this particular choice of objective function, we have that $\nabla f(x_{k-1})  = - (x^\star - x_{k-1} )$. Hence, the update in Eq.\ (\ref{eq:FWupdate}) becomes
\[
p_k \in \arg \min_{u \in \mathcal P } -(x^\star -x_{k-1} )^T u,
\] 
which is equivalent to solving the linear program \texttt{LP}  $( -(x^\star- x_{k-1}), \mathcal B)$. That is, computing an extreme point with {Frank-Wolfe} and \texttt{Birkhoff} is the same except for the ceiling.\footnote{Recall also that with FW there is not requirement that $x^\star \succeq x$.} Note that by ceiling the vector $-(x^\star- x_{k-1})$, we are ``weighting'' all the components that are not equal to zero equally. Without the ceiling, the Frank-Wolfe update takes into account the geometry of the decomposition, i.e., how close $x_{k-1}$ is to $x^\star$ entry-wise.
\label{ob:1}
\end{observation}

\begin{observation}[Steepest descent permutation] The extreme points selected by Frank-Wolfe corresponds to obtaining the ``steepest'' descent direction, or direction $(p_k-x_{k-1})$ that has the smallest angle with respect to $(x^\star -x_{k-1})$. Note that $(x^\star -x_{k-1}) = - \nabla f(x_{k-1})$ is the direction that goes straight to the target value $x^\star$, and that
\begin{align*}
& \arg \min_{u \in \mathcal P } \nabla f(x_{k-1})^T u \\
& \quad \stackrel{\text{(a)}}{=}  \arg \min_{u \in \mathcal P } \| \nabla f(x_{k-1}) \|_2 \| u \|_2 \cos \phi_{\langle \nabla f , u \rangle} \\
& \quad \stackrel{\text{(b)}}{=} \arg \min_{u \in \mathcal P } \ \cos \phi_{\langle \nabla f , u \rangle} 
\end{align*}
where (a) follows from the dot product and (b) since $\| p \|_2 = \sqrt n$ for all $p \in \mathcal P$, and $ \| \nabla f(x_{k-1}) \|_2$ does not depend on $p$. The RHS of the last equation corresponds to maximizing $\cos \phi_{\langle - \nabla f , p \rangle} $, which is equivalent to finding the $p \in \mathcal P$ that minimizes the angle between $- \nabla f(x) = (x^\star - x)$ and $(p-x)$. 
Furthermore, since the Birkhoff polytope is regular and the number of extreme points increases factorially with $n$, we can expect $ \phi_{\langle -\nabla f , p \rangle}$ to be small. Figure \ref{fig:steepest_descent} shows, schematically, how Frank-Wolfe selects the extreme point that has the smallest angle with respect to $(x^\star - x)$. The black dots with a red cross are ``non-descent'' permutations that will not improve the decomposition approximation. 
\label{ob:2}
\end{observation}

Both observations rely on the objective function being quadratic; however, we can expect similar properties for other smooth convex objectives. For example, we could use $f(x) = (x^\star-x)^T Q (x^\star-x) $ where $Q$ is a positive semi-definite matrix that emphasizes which of the components in vector $x^\star-x$ to minimize. In Section \ref{sec:birkhoff+}, we will include a $\log$-barrier function to the objective. In the rest of the section, we will use a quadratic objective function to streamline exposition.

\begin{algorithm}[t]
\caption{Frank-Wolfe (\texttt{FW}) with quadratic objective and line search}
\begin{algorithmic}[1]
\STATE As Algorithm \ref{al:birkhoff_algorithm_fwform}, but set $x_0 \in \mathcal P$ and replace lines $\circ$, {\small \ding{73}}, $\diamond$ with 
\STATE   $\quad p_k \leftarrow$ \texttt{LP}  $(-(x^\star - x_{k-1}) , \mathcal B)$\\ 
\STATE   $\quad \theta_k \leftarrow  (x^\star-x_{k-1})^T (p_k-x_{k-1}) / \| p_k -x_{k-1} \|_2^2   $ \label{lst:FWstep}
\STATE   $\quad x_k \leftarrow  x_{k-1} + \theta_k (p_k-x_{k-1})$
\end{algorithmic}
\label{al:FW}
\end{algorithm}

\subsection{Frank-Wolfe with line search}
\label{sec:FW_linesearch}
The procedure of the Frank-Wolfe algorithm is given in Algorithm \ref{al:FW}. Differently from Birkhoff's approach, \texttt{FW} uses an extreme point as a starting point instead of the origin. Note that $0 \notin \mathcal B$. The choice of step size is indicated in step 3, and corresponds to carrying out line search. This can be easily verified.  Let $x_k \coloneqq x_{k-1} + \theta_k (p_k - x_{k-1})$ be the $k$'th iterate, and observe that we can write 
\begin{align*}
& \frac{1}{2} \| x_{k} - x^\star \|_2^2 - \frac{1}{2} \| x_{k-1} - x^\star \|_2^2\\
&  = \frac{1}{2} \| x_{k-1} + \theta_k (p_k - x_{k-1}) - x^\star \|_2^2 - \frac{1}{2} \| x_{k-1} - x^\star \|_2^2 \\
& = \theta_k (x_{k-1} - x^\star)^T (p_k - x_{k-1}) + \frac{\theta^2}{2} \| p_k - x_{k-1} \|_2^2.
\end{align*}
The RHS of the last equation is a quadratic function in $\theta_k$, whose minimizer can be obtained in closed form. And since equality holds in the last equation, minimizing the quadratic function on the RHS is equivalent to minimizing the LHS with line search. 
Hence, Algorithm \ref{al:FW} corresponds to Frank-Wolfe with line search, and so from  \cite[Theorem 1]{Jag13},\footnote{The bound in Eq.\ \eqref{eq:FWconvergence} follows from Theorem 1 in \cite{Jag13} with $\delta = 0$ (i.e., in our problem the gradients are noiseless) and $C_f = \max_{u,v \in \mathcal B} \| u-v\|_2^2 = 2n$.} we have the bound
\begin{align}
\| x_k - x^\star \|_2^2 \le \frac{4n}{k+2}. \label{eq:FWconvergence}
\end{align}
By rearranging terms in Eq.\ \eqref{eq:FWconvergence}, we can obtain an upper bound on the sparsity of \texttt{FW}.
\begin{theorem}[\texttt{FW} sparisty]
\label{th:FW}
Algorithm \ref{al:FW} obtains an $\epsilon$-approximate decomposition with at most $k \le 4n / \epsilon^2$ permutation matrices, where $\epsilon  = \| x_k - x^\star \|_2$.
\end{theorem}
The bound in Theorem \ref{th:FW} says that the sparsity increases exponentially with the error, and so it does not allow us to obtain good approximations that are also sparse. 
One of the issues with first-order-methods is the zig-zagging phenomenon\footnote{See the discussion on page 2 in \cite{LJ15}.} when the approximate decomposition is close to $x^\star$.  Hence, even though \texttt{FW} selects the steepest descent direction, the choice of step size is not enough. 
One way to avoid zig-zagging is to recompute the weights of all the atoms or extreme points discovered so far, which is in essence what the fully corrective variant of the algorithm does. 

 \begin{algorithm}[t]
\caption{Fully Corrective Frank-Wolfe (\texttt{FCFW})}
\begin{algorithmic}[1]
\STATE As Algorithm \ref{al:birkhoff_algorithm_fwform}, but set $x_0 \in \mathcal P$ and define  $V_0 = \emptyset$. Let $\Delta_k$ be the $k$-simplex. Replace lines $\circ$, {\small \ding{73}}, $\diamond$ with 
\STATE   $\quad p_k \leftarrow$ \texttt{LP}  $(\nabla f(x_{k-1}) , \mathcal B)$\\ 
\STATE   $\quad V_k \leftarrow [ V_{k-1} , p_k ]$
\STATE   $\quad (\theta_1,\dots,\theta_k) \leftarrow  \arg \min_{u \in \Delta_k}\ \| V_k u - x^\star \|^2_2 $ 
\STATE   $\quad x_k \leftarrow V_k (\theta_1,\dots,\theta_k)$
\end{algorithmic}
\label{al:FCFW}
\end{algorithm}

\subsection{Fully Corrective Frank-Wolfe (FCFW)}
\label{sec:FCFW}
This variant of Frank-Wolfe differers from the classic algorithm because it provides the best approximation with the number of extreme points selected up to iteration $k$. The \texttt{FCFW} procedure is described in Algorithm \ref{al:FCFW}. As in the \texttt{FW} algorithm, it starts from an arbitrary $p \in \mathcal P$ and computes a new permutation by solving a linear program \texttt{LP}  $(\nabla f(x_{k-1}) , \mathcal B)$. The main difference is that the permutations are collected in matrix $V_k$, and the weights $(\theta_1,\dots,\theta_k)$ selected to minimize $\| V_k (\theta_1,\dots,\theta_k) - x^\star \|_2^2$ subject to $\sum_{i=1}^k \theta_ i = 1$ and $\theta_i \ge 0$ for all $i=1,\dots,k$. 
An important difference of Algorithm \ref{al:FCFW} with respect to Algorithm \ref{al:FW} is that the computation of a new collection of weights involves solving a quadratic program (QP) whose dimension increases with the number of iterations. 
For this Frank-Wolfe variant, from Theorem 1 in \cite{LJ15}, we have the bound
\begin{align}
 \|  x_k - x^\star \|_2^2 \le \| x_0 - x^\star \|_2^2 \exp \left( - \frac{\mu}{4L} \left(\frac{\lambda}{M}\right)^2 k \right), \label{eq:FCFWconvergence}
\end{align}
where $\mu/L$ is the condition number and $(\lambda/M)^2$ the eccentricity\footnote{The eccentricity of a set is  similar to the condition number of a function; see \cite[pp. 461]{BV04}} of the Birkhoff polytope. These two parameters are usually not known, however, not in our problem since the Birkhoff polytope and the objective function $f(x) = (1/2) \| x - x^\star  \|_2^2$ have remarkable structure. We establish the eccentricity of the Birkhoff polytope in the next lemma. 

\begin{lemma}
\label{th:eccentricity}
The eccentricity $(\lambda/M)^2$ of the Birkhoff polytope is lower bounded by ${1}/(2n^3)$.  
\end{lemma}
Using the last lemma and the fact that the condition number  of the objective function ($\mu /L$) is equal to $1$, we can obtain the \texttt{FCFW}'s sparsity.
\begin{theorem}[\texttt{FCFW} sparsity]
Algorithm \ref{al:FCFW} obtains an $\epsilon$-approximate decomposition with at most $k \le 8n^3 \log(2n / \epsilon^2)$ permutation matrices, where $\epsilon = \| x_k - x^\star \|_2$.
\label{th:FCFW}
\end{theorem}

From the last theorem, we have that the number of extreme points required to obtain an approximate Birkhoff decomposition increases logarithmically with the decomposition error. This is a huge improvement with respect to the sparsity result obtained with the line search \texttt{FW} in Theorem \ref{th:FW}.
Unfortunately, \texttt{FCFW} is less exciting in practice because recomputing the weights is expensive computationally since the size of the quadratic program (step 4 in Algorithm \ref{al:FCFW}) increases with the number of permutations. Furthermore, the accuracies of the quadratic solvers such as SCS~\cite{OCP+16}, Ipopt \cite{Ipopt}, and Gurobi~\cite{Gurobi} are in the order of $10^{-6}$, which means that we cannot obtain decompositions with accuracies below $10^{-3}$. The latter can be observed in Figure \ref{fig:sparsity}a in the numerical evaluation. 


\section{New Algorithm}
\label{sec:birkhoff+}

\texttt{Birkhoff} and \texttt{FCFW} algorithms have both logarithmic sparsity, but they are very different algorithmically. On the one hand, weights are easy to compute in Birkhoff's approach,\footnote{Birkhoff's step size requires to find the smallest of $n$ elements.} but finding a good permutation matrix is slow as it requires to solve \emph{multiple} linear programs (e.g., \cite{LML+15}). In contrast, \texttt{FCFW} can obtain a good permutation matrix by solving a \emph{single} linear program (see Observation \ref{ob:2}), but it requires to solve a quadratic program to (re)calculate the weights.

In this section, we present \texttt{Birkhoff+} (Algorithm \ref{al:birkhoff+}), a variation of the original Birkhoff's algorithm that uses the intuition behind Frank-Wolfe to obtain sparse decompositions in a fast manner. The performance of \texttt{Birkhoff+}  is evaluated in Section \ref{sec:numerical_evaluation}.

 \begin{figure}
\centering
\includegraphics[width=0.49\textwidth]{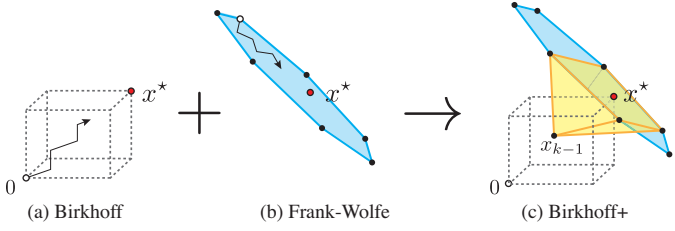}
\caption{Schematic illustration of the (a) Birkhoff's and (b) Frank-Wolfe approaches and how they combine into the (c) new setup. The yellow polygon represents the convex hull of $x_{k-1}$ and the permutation matrices in $\mathcal I_k(\alpha)$.}
\label{fig:cones}
\end{figure}

\subsection{Approach}
\label{sec:approach}

The intuition behind our approach is shown schematically in Figure \ref{fig:cones}. In brief, Birkhoff's algorithm (Figure \ref{fig:cones}a) can be seen as constructing a path from the origin ($x_0 = 0$) to the target value ($x^\star$) while always remaining in the dotted box (i.e., $x_k \preceq x^\star$ for all $k=0,1,2,\dots$). Frank-Wolfe (Figure \ref{fig:cones}b), on the other hand, constructs a path from a permutation matrix $x_0 \in \mathcal P$ to the target value $x^\star$ within the polytope of doubly stochastic matrices (blue surface). Our approach (Figure \ref{fig:cones}c) can be regarded as using Frank-Wolfe within the polytope $\conv{(\mathcal I_k(\alpha) \cup x_{k-1}})$ (yellow polygon in Figure \ref{fig:cones}c) with the additional constraint that the approximate decomposition must be within the dotted box. That is, we want to use the path or permutations that Frank-Wolfe would select while remaining in the box that characterizes the Birkhoff's approach. It is important to use $\conv{(\mathcal I_k(\alpha) \cup x_{k-1}})$ instead of $\conv{( \mathcal P \cup x_{k-1}})$ (i.e., all permutations) as  the algorithm may otherwise not converge. The latter is shown formally in the following theorem.

\begin{theorem}
\label{th:birkhoff_algorithm_fwform_notcvg}
Consider Algorithm \ref{al:birkhoff_algorithm_fwform} and replace line $\circ$ with $\texttt{LP} (-(x^\star - x_{k-1}))$. Then, there may not exist a $k$ for which $\| x_k - x^\star \|_2^2 \le \epsilon$ for any $\epsilon >0$.
\end{theorem}
We can prove the theorem by example. Suppose we want to decompose the following $n \times n$ doubly stochastic matrix 
\begin{align}
\begin{bmatrix}
1- 1/n & 0& \cdots  & 0 &1/n \\
0 & 1 -1/n &  &  & 1/n\\
\vdots & &   \ddots && \vdots \\
 0& &   & 1 - 1/n&  1/n\\
 1/n & 1/n & \cdots & 1/n & 0 \\
\end{bmatrix}
\label{eq:matrix_dd}
\end{align}
That is,  (i) the first $n-1$ entries in the diagonal are equal to $1 - 1/n$, (ii) the first $n-1$ entries of the last row are equal to $1/n$, and (iii) the first $n-1$ entries of the last column are equal to $1/n$. 
Note the sum of each row and column is equal to one. Next, suppose that $f(x) = (1/2) \|x^\star -x \|_2^2$ where $x^\star$ is the matrix in Eq.\ \eqref{eq:matrix_dd} in column form.  In the first iteration ($x_0= 0$), Frank-Wolfe selects a permutation by solving the linear program $\texttt{LP}(-x^\star, \mathcal B)$, the solution of which is the identity matrix since the doubly stochastic matrix in Eq.\ \eqref{eq:matrix_dd} is diagonally dominant.  And because the last entry of the matrix in Eq.\ \eqref{eq:matrix_dd} is equal to zero, we  have that $\texttt{BIRKHOFF\textunderscore STEP} (x^\star, x_{k-1}, p_k) = 0$ and therefore $x_k = x_{k-1}$. That is, the algorithm will be  ``stuck.'' 

In sum, a Birkhoff-type algorithm that selects permutation matrices with Frank-Wolfe using all the permutation matrices $\mathcal P$ may not converge.  However, we can use Frank-Wolfe with the permutations in the set  $\mathcal I_k(\alpha)$, which ensures not only that the algorithm converges but that this has logarithmic sparsity (Theorem \ref{th:birkhoff_linear_convergence}).

\subsubsection{Objective function with barrier}

Since Birkhoff's approach restricts $x_{k}$ to remain in the Birkhoff's dotted box (see Figure \ref{fig:cones}), it is reasonable to use an objective function that aims to construct a path to $x^\star$ from within the box. For that, we define
\begin{align}
f_\beta (x) = f(x) - \beta \sum_{j=1}^d \log(x^\star(j) - x(j) + \epsilon / d ), \label{eq:fbeta}
\end{align}
where $\beta \ge 0$ and $x(j)$ is the $j$'th component of vector $x$. Note that $f_\beta$ is convex as this is the composition of $f$ plus a convex penalty/barrier function $-  \beta \sum_{j=1}\log(x^\star(j) - x_{k-1} (j) + \epsilon / d )$. The term $\epsilon / d$ in the barrier is used for numerical stability as otherwise the barrier goes to $+\infty$ when $x^\star (j) = x(j)$. The motivation for using a barrier function comes from interior point methods in optimization, where parameter $\beta$ is typically tuned throughout  the algorithm to allow $x_k \to x^\star$. Note that $f_\beta \to f$ as $\beta \to 0$.

\subsection{\texttt{Birkhoff+} algorithm description and complexity}

The procedure of \texttt{Birkhoff+} is described in Algorithm \ref{al:birkhoff+}, and consists of replacing how permutation matrices are selected in Algorithm \ref{al:birkhoff_algorithm_fwform} with $\texttt{LP}(\nabla f_\beta (x_{k-1}), \conv{(\mathcal I_k(\alpha) )})$, where   $f_\beta$ is as defined in Eq.\ \eqref{eq:fbeta}.  Parameter $\beta$ can be selected to emphasize the barrier over the objective function $f$. In our case, we do not need $\beta \to 0$ as by selecting permutations from $\mathcal I_k(\alpha)$ is enough to allow the algorithm to make progress. The convergence of the algorithm is stated formally in the following corollary.

\begin{corollary}
\label{th:bir_coro}
Algorithm \ref{al:birkhoff+} obtains an $\epsilon$-approximate decomposition with at most $k \le O(\log(1 / \epsilon ))$ permutation matrices. 
\end{corollary}

The complexity of \texttt{Birkhoff+} per iteration is equal to solving a linear program with a simplex type method. The linear program $\texttt{LP}(\nabla f_\beta (x_{k-1}), \conv{(\mathcal I_k(\alpha)})$ can be carried out with $\texttt{LP}(\nabla f_\beta (x_{k-1}) + b_k, \mathcal B)$ where $b_k  = d/\epsilon \cdot \mathbb I_{\{0,1\}}( x^\star - x_{k-1} \preceq  \alpha)$ is a penalty vector to force the solver to do not select the components of vector $(x^\star - x_k)$ smaller than $\alpha$.

Finally, we note that \texttt{Birkhoff+} depends on how we define set $\mathcal I_k(\alpha)$. Algorithm \ref{al:birkhoff_refinement} is a meta-heuristic for selecting $\alpha$ based on Birkhoff's step size. 
In particular, $\alpha$ is set to $(1-\sum_{i=1}^k \theta_i)/n^2$ in the first iteration and then equal to the largest step size for the permutation selected using the Frank-Wolfe-type update. The search for a large $\alpha$ terminates when the maximum number of repetitions $(\texttt{max\_rep})$ is reached or the value of $\alpha$ does not increase. We call Algorithm \ref{al:birkhoff_refinement} \texttt{Birkhoff+(\#)}, where \texttt{\#} indicates the maximum number of permutation refinements. \texttt{Birkhoff+(1)} is equivalent to \texttt{Birkhoff+} as it computes only one permutation matrix. 

\begin{algorithm}[t]
\caption{\texttt{Birkhoff+} }
\begin{algorithmic}[1]
\STATE As Algorithm \ref{al:birkhoff_algorithm_fwform}, but take $\beta \ge 0$ also as input. Replace line $\circ$ with
\STATE \quad $\alpha \leftarrow (1-\sum_{i=1}^{k-1} \theta_i ) / n^2 $
\STATE \quad $p_k \leftarrow \texttt{LP}(\nabla f_\beta (x_{k-1}), \conv{(\mathcal I_k(\alpha))})$ 
\end{algorithmic}
\label{al:birkhoff+}
\end{algorithm}

\begin{algorithm}[t]
\caption{\texttt{Birkhoff+(max\textunderscore rep)} --- with permutation selection refinement }
\begin{algorithmic}[1]
\STATE As Algorithm \ref{al:birkhoff+}, but replace line $\circ$ with
\FOR {$i =1, \dots, \texttt{max\textunderscore rep}$} 
	\STATE $p_i \leftarrow \texttt{LP}(\nabla f_\beta (x_{k-1}), \conv{(\mathcal I_k(\alpha)})$
	\STATE $\theta_i \leftarrow \texttt{BIRKHOFF\textunderscore STEP} (x^\star, x_{k-1}, p_k)$  
	\STATE  \textbf{if} ($\theta_i > \alpha$) $\alpha \leftarrow \texttt{BIRKHOFF\textunderscore STEP} (x^\star, x_{k-1}, p_k)$
	\STATE  \textbf{else} exit while loop
	\STATE $p_k \leftarrow p_i$
 \ENDFOR
  \end{algorithmic}
\label{al:birkhoff_refinement}
\end{algorithm}


\section{Numerical Evaluation} 
\label{sec:numerical_evaluation}

In this section, we evaluate performance of \texttt{Birkhoff+} and compare it to existing algorithms. Our goal is to illustrate the algorithms' characteristics and how different traffic matrices affect the performance of a circuit switch in terms of throughput, configurations computation time, and number of configurations. The code of \texttt{Birkhoff+} is available as a Julia \cite{Julia} package in \cite{BirkhoffDecompositionPackage}.

\subsection{Setup}
\label{sec:evaluation_setup}

The  \texttt{Birkhoff}, \texttt{FW}, \texttt{FCFW},  \texttt{Birkhoff+} and \texttt{Birkhoff+(\#)} algorithms are implemented in Julia \cite{Julia} and as indicated in Algorithms \ref{al:birkhoff_algorithm_fwform}, \ref{al:FW}, \ref{al:FCFW}, \ref{al:birkhoff+} and \ref{al:birkhoff_refinement} respectively. Parameter $\beta$ is fixed to $1$ and the maximum number of permutation refinements in \texttt{Birkhoff+(\#)}  to $10$ --- however, we observe in the experiments that the actual number of permutation refinements is usually less than $3$. \texttt{Solstice} corresponds to Algorithm 2 in \cite{LML+15}, and \texttt{Eclipse} to Algorithm 2 in \cite{BAV16}. The linear programs $\texttt{LP}(\cdot ,\cdot)$  are carried out with Clp \cite{Clp} in all algorithms and return an extreme point/permutation matrix. The quadratic programs in the \texttt{FCFW} algorithm are carried out with Ipopt \cite{Ipopt}. Both solvers are open-source. 

Traffic demand matrices are generated by sampling permutations uniformly at random, and weights are selected to model the type of load in data centers. In particular, we follow the evaluation scenario in  \cite{BAV16}, where traffic matrices are sparse and consist of 12 flows. Three of the flows are large and carry the 70\% of the load, while the rest are small flows and carry the remaining 30\% of the traffic. 
We note that the traffic matrix in practical scenarios may be below the switch’s capacity (i.e., the sum of each row or column may be smaller than one), and so we first need to add a virtual load to the traffic matrix to make it doubly stochastic.\footnote{The work in \cite{LML+15} (see Section 4.2.1) uses the term ``stuffing'' for adding virtual load to the traffic matrix. Stuffing can be seen as a special type of projection of the demand matrix onto the Birkhoff polytope. Technically, for a demand matrix $D$, we need to find a matrix $S \in \mathcal B - D$. Matrix $S$ may not be unique and finding the best virtual load matrix for our algorithm is an interesting problem but out of the scope of the paper.} For simplicity, we assume in the evaluation that the demand matrices are doubly stochastic. 

Finally, the numerical evaluation is carried out on a computer equipped with an Intel i7 8700B (3.2 GHz) CPU and 32 GB of memory. The version of Julia is 1.3.1.

\subsection{Experiments}

We first study the algorithms' characteristics, and then show how those affect the performance of a circuit switch.

\subsubsection{Decomposition approximation vs. number of permutations and time} 
\label{sec:zero_exp}

\begin{figure}
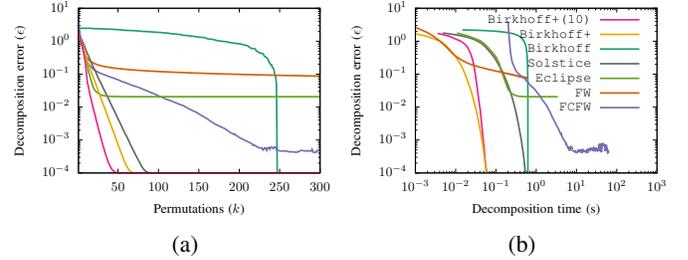

\begin{tabular}{cc}
\hspace{-1em} {\resizebox{0.51\columnwidth}{!}{\input{Figures/EXP0_RC_32.tex}}} & \hspace{-1em}
\hspace{-1em} {\resizebox{0.51\columnwidth}{!}{\input{Figures/EXP0_RC_32_speed.tex}}} \\
\small (a) & \hspace{-1em} \small (b) 
\end{tabular}
\caption{Decomposition error ($\epsilon$) of \texttt{Birkhoff}, \texttt{FW}, \texttt{FCFW}, \texttt{Solstice}, \texttt{Eclipse}, \texttt{Birkhoff+}, and \texttt{Birkhoff+(10)} algorithms depending on (a) the number or permutations and (b) time. The figure shows the average of 50 realizations.   }
\label{fig:sparsity}
\end{figure}

We set $n = 32$ and sample traffic matrices as indicated in Section \ref{sec:evaluation_setup}. Also, we fix $\epsilon = 10^{-4}$, $k_\text{max} = 300$ and $\delta = 10^{-2}$ (just for \texttt{Eclipse})\footnote{The value corresponds to having a switching cost of 10 ms.}.  Figure \ref{fig:sparsity} shows the algorithms decomposition error in terms of permutations and time. The results are the average of 50 realizations. 

Observe from Figure \ref{fig:sparsity}a that the decomposition error of \texttt{Birkhoff} is large until it converges exactly in the last iteration ($k \approx 250$). On the other hand, \texttt{FW} progresses quickly, but  it slows down drastically around $\epsilon = 0.9$. The latter is due to the $O(1/\epsilon^2)$ sparsity rate and the zig-zagging phenomenon typical in first-order-methods (see Section \ref{sec:FW_linesearch}). The \texttt{FCFW} has a better sparsity performance than \texttt{FW}, but it cannot obtain decomposition with an $\epsilon$ below $0.5 \cdot 10^{-3}$ due to the numerical accuracy of the quadratic solvers (see Section \ref{sec:FCFW}). \texttt{Eclipse} has a better performance than previous algorithms until it gets stuck between $\epsilon \in [10^{-2}, 10^{-1}  ]$. We conjecture the latter is because \texttt{Eclipse} selects permutations using a Max-Weight-type matching, and so it may face similar issues as when we combine Frank-Wolfe and Birkhoff approaches directly; see discussion in Section \ref{sec:approach}. Also, the performance guarantees of \texttt{Eclipse} given in 
\cite{BAV16}  are for the problem type in \cite{LH03} (see Section \ref{sec:related_work}) and not for decomposing a doubly stochastic matrix. 
Finally, observe that \texttt{Solstice}, \texttt{Birkhoff+}, and \texttt{Birkhoff+(10)\footnote{The number in the parentheses is the maximum number of permutations refinements (\texttt{max\textunderscore rep}).}} have all better sparsity performance than the previous algorithms and that $\texttt{Birkhoff+(10)}$ is noticeably better for $\epsilon < 0.1$. The last three algorithms have linear convergence/logarithmic sparsity ($y$-axis is in log-scale) but different condition numbers:  $0.89$, $0.85$ and $0.82$ respectively.\footnote{Average of the $50$ first iterations.} Recall the condition number indicates how an additional permutation reduces the decomposition error multiplicatively (see discussion in Section \ref{sec:convergence}).

Figure \ref{fig:sparsity}b shows the decomposition error against the running time. Observe that \texttt{Birkhoff+} is the fastest followed by \texttt{Birkhoff+(10)}. \texttt{FW} is also fast for $\epsilon > 0.1$, but it slows down afterward for the same reason explained above. \texttt{Solstice} and \texttt{Eclipse} are both slower than \texttt{Birkhoff+} by an order of magnitude since they need to solve multiple linear programs to select a permutation matrix. The running time of \texttt{Birkhoff} is in line with the sparsity results: it makes slow progress until it converges exactly in the last iteration. Finally, \texttt{FCFW} is the slowest as it has to recompute all the weights (i.e., solve a quadratic program) every time it adds a new permutation to the decomposition.

\begin{figure*}[t!]
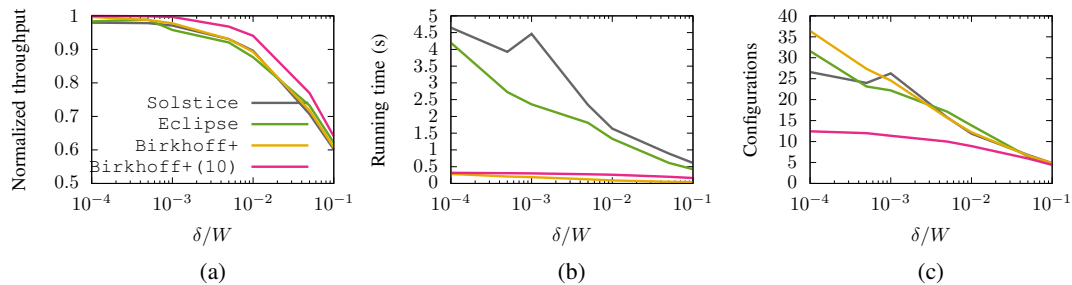

\centering
\begin{tabular}{ccc}
{\resizebox{0.24\textwidth}{!}{\input{Figures/RC_100.tex}}} & 
{\resizebox{0.24\textwidth}{!}{\input{Figures/RT_100.tex}}} & 
{\resizebox{0.24\textwidth}{!}{\input{Figures/C_100.tex}}} \\
\small (a) & \small (b) & \small (c)
\end{tabular}
\caption{Circuit switch performance (throughput, running time, and number of configurations) depending on $\delta/ W$, where $\delta$ is the switching time and $W$ the time window duration. The figures show the average of 50 realizations.}
\label{fig:ratio}
\end{figure*}

\subsubsection{Circuit switch performance}  We now evaluate the algorithm's performance when used to compute the switching configurations for a circuit switch with $n = 100$ ports. The performance metrics we evaluate are the throughput, the configurations computation time, and the number of switching configurations. We carry out three experiments where we vary the reconfiguration cost, the skewness and sparsity of the traffic matrix, and the configurations computation overhead. 
Importantly, now the traffic matrix $X^\star$ is associated with a time window $W$ that enforces the decomposition to satisfy $\sum_{i=1}^k (\theta_i + \delta) \le W$, i.e., the time spent transmitting $(\sum_{i=1}^k \theta_i)$ and reconfiguring ($\delta k$) cannot exceed the time window duration ($W$).  
Finally, we only evaluate \texttt{Solstice}, \texttt{Eclipse}, \texttt{Birkhoff+}, and \texttt{Birkhoff+(10)} as (i) \texttt{Birkhoff} and \texttt{FW} have a poor performance, and (ii) \texttt{FCFW} is very slow when $n \ge 32$ (see times in Figure\ \ref{fig:sparsity}b).

\textbf{Experiment 1 (impact of reconfiguration time).}  
Figure \ref{fig:ratio} shows the algorithms' performance in terms of throughput, running time, and the number of configurations depending on the ratio $\delta / W$ (the impact of the reconfiguration delay proportionally to the time window duration). Observe from the figure that \texttt{Birkhoff+(10)} outperforms the other algorithms in terms of throughput. For instance, for $\delta/W = 10^{-2}$, \texttt{Birkhoff+(10)} achieves a $7\%$ more throughput than \texttt{Eclipse} and \texttt{Solstice}. \texttt{Birkhoff+} has almost the same throughput than  \texttt{Eclipse} and \texttt{Solstice}. Regarding the time required to compute the switching configurations,  \texttt{Solstice} and \texttt{Eclipse} are slower than \texttt{Birkhoff+} and \texttt{Birkhoff+(10)} by an order of magnitude; however, the difference decreases as $\delta / W$ increases because we have fewer switching configurations as a result of larger reconfiguration penalties (c.f.\  Figure \ref{fig:ratio}b and Figure \ref{fig:ratio}c). Finally, observe from Figure \ref{fig:ratio}c that \texttt{Birkhoff+(10)} can obtain decompositions with half of the configurations compared to other algorithms when $\epsilon$ is small (i.e., $10^{-4}$). \textbf{Conclusions:} \texttt{Birkhoff+} has the same performance in terms of throughput and number of switching configurations than \texttt{Solstice} and \texttt{Eclipse}, but it is 10 times faster. \texttt{Birkhoff+(10)} obtains higher throughput than all algorithms and it is only slightly slower than \texttt{Birkhoff+}.

\begin{figure*}
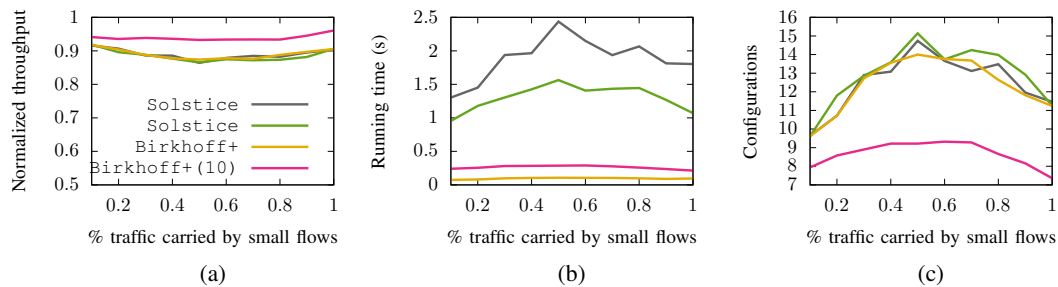

\centering
\begin{tabular}{ccc}
{\resizebox{0.24\textwidth}{!}{\input{Figures/EXP2_RC_100.tex}}} & 
{\resizebox{0.24\textwidth}{!}{\input{Figures/EXP2_RT_100.tex}}} & 
{\resizebox{0.24\textwidth}{!}{\input{Figures/EXP2_C_100.tex}}} \\
\small (a) & \small (b) & \small (c) 
\end{tabular}
\caption{Circuit switch performance (throughput, running time, and number of configurations) depending on the load carried by the small flows. The figures show the average of 50 realizations.}\label{fig:skew}
\vspace{-1em}
\end{figure*}

\textbf{Experiment 2 (sparsity and skewness).}  
Now we set $\delta / W = 10^{-2}$ and evaluate the algorithms' performance depending on the skewness and the sparsity of the demand matrix. In Figure \ref{fig:skew}, we vary the fraction of the load   carried by the small flows. Observe that as before, \texttt{Birkhoff+(10)} outperforms the other algorithms, and that \texttt{Birkhoff+},  \texttt{Solstice}, and \texttt{Eclipse} are almost the same in terms of throughput for different demand matrices. Furthermore, there is little variation on the running time and the number of switching configurations---despite a slight bend in the curves when the traffic matrix contains the same fraction of large and small flows. 

\begin{figure*}
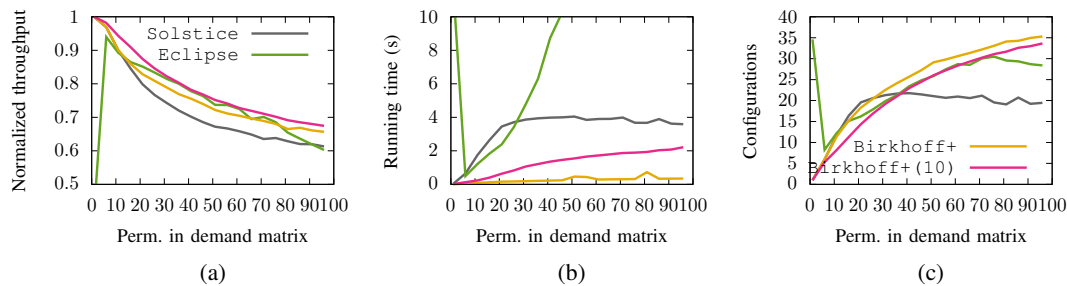

\centering
\begin{tabular}{ccc}
{\resizebox{0.24\textwidth}{!}{\input{Figures/EXP2B_RC_100.tex}}} & 
{\resizebox{0.24\textwidth}{!}{\input{Figures/EXP2B_RT_100.tex}}} & 
{\resizebox{0.24\textwidth}{!}{\input{Figures/EXP2B_C_100.tex}}} \\
\small (a) & \small (b) & \small (c)
\end{tabular}
\caption{Circuit switch performance (throughput, running time, and number of configurations) depending on the number of permutations matrices used to generate the traffic matrix. The figures shows the average of 50 realizations.}
\label{fig:sparse}
\end{figure*}

In Figure \ref{fig:sparse}, we show the results when we vary the number of permutations used to generate the demand matrix. Each permutation matrix is sampled as explained in Section \ref{sec:evaluation_setup}. Observe from the figure that the sparsity of the traffic matrix has a significant impact on the throughput, running time, and the number of the switching configurations. In particular, observe from Figure \ref{fig:sparse}a that the throughput of all algorithms decreases and that \texttt{Eclipse} is comparable to \texttt{Birkhoff+(10)} as the traffic matrix becomes denser. However, the running time of \texttt{Birkhoff+(10)} does not explode (see Figure \ref{fig:sparse}b) and \texttt{Birkhoff+(10)} does not get stuck when the traffic demand matrix is very sparse.\footnote{See discussion in Section \ref{sec:zero_exp}.} 
Regarding \texttt{Birkhoff+}, observe that now the running times difference with \texttt{Birkhoff+(10)} becomes more noticeable as the demand matrix becomes denser. Finally, observe from Figure \ref{fig:sparse}c that the number of switching configurations increases with the density of the traffic matrix for all algorithms. 
\textbf{Conclusions:} The skewness of the demand matrix has little impact on to the performance of all algorithms. The sparsity, on the other hand, plays an important role. \texttt{Eclipse} has a similar performance than \texttt{Birkhoff+(10)}, but it is notably slower. 


\begin{figure}
\centering
{\resizebox{0.45\columnwidth}{!}{\input{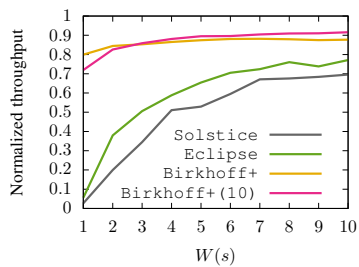}}}
\caption{Circuit switch throughput when the time to compute the switching configurations is an overhead. The figure shows the average of 50 realizations.}\label{fig:txtx}
\end{figure}

\textbf{Experiment 3 (configurations computation overhead).} 
This experiment shows how time to compute the switching configurations affects the circuit switch's throughput. In particular, we set $\delta / W = 10^{-2}$ and truncate the decomposition to satisfy $\sum_{i=1}^k (\theta_i + \delta) \le W - T$, where $T$ is the time to compute the switching configurations. The values of $T$ for this particular setting are given in Figure \ref{fig:ratio}b. Figure \ref{fig:txtx} shows the throughput for different values of $W$ in seconds. Observe from the figure that the throughput increases with $W$ for all algorithms. When $W$ is small (i.e., the decomposition computation overhead is large), \texttt{Birkhoff+} has a higher throughput than \texttt{Birkhoff+(10)} because it is faster---recall \texttt{Birkhoff+} selects a new switching configuration by solving a single linear program. However, \texttt{Birkhoff+(10)}'s throughput is higher when $W > 5$ since the reconfiguration time is larger than the time to compute the  switching configurations. 
Regarding \texttt{Eclipse} and \texttt{Solstice}, observe that both are affected heavily by the decomposition overhead. For instance, when $W = 5$, \texttt{Birkhoff+} has 67\% and 34\% more throughput than \texttt{Solstice} and \texttt{Eclipse} respectively. Also, note that when $W= 1$,  \texttt{Birkhoff+} can serve 80\% of the traffic whereas \texttt{Solstice} and \texttt{Eclipse} almost nothing. 
\textbf{Conclusions:} \texttt{Birkhoff+} outperforms  \texttt{Solstice} and \texttt{Eclipse} and it is slightly better than  \texttt{Birkhoff+(10)} when the time windows are short. As with the reconfiguration costs, the benefit of computing switching configurations fast diminishes as the time window duration increases. 

\section{Conclusions}
\label{sec:conclusions}

This paper studies how to compute switching configurations for circuit switches. We have revisited Birkhoff's approach and established its properties in terms of the number of switching configurations required to obtain an approximate representation of a traffic matrix. A new algorithm (\texttt{Birkhoff+}) is proposed, which obtains representations with fewer switching configurations than previous work (\texttt{Solstice}, \texttt{Eclipse}) and is 10-100 times faster depending on the setting. The latter is important in terms of throughput when traffic bursts are short-lived, and so the time required to compute the switching configurations is a non-negligible overhead. We also propose a variant of \texttt{Birkhoff+} that is slightly slower but obtains representations with even fewer switching configurations. The performance of the proposed algorithms is evaluated through exhaustive numerical experiments for traffic demand matrices that capture traffic characteristics in data centers.

\section{Acknowledgements}

This work has received funding from the European Union’s Horizon 2020 research and innovation programme under the Marie Skłodowska-Curie grant agreement No.\ 795244. 

This work was sponsored by the NSF Award 1815676. This research was supported by ARO (W911NF1810378).

The first author would like to thank Ehsan Kazemi (Yale Institute for Network Science) for many helpful conversations and the counter-example in the proof of Theorem \ref{th:birkhoff_algorithm_fwform_notcvg}. 


\bibliographystyle{IEEEtran}
\bibliography{references}

\appendix
\section{Appendix}
\label{sec:proofs}

\subsection{Proofs of Section \ref{sec:birkhoff_revisited}}

We start by presenting two lemmas. The first lemma gives an upper bound on the Frobenius norm of a doubly stochastic matrix.
\begin{lemma}
$\| X \|_F \le \sqrt n$ for any doubly stochastic matrix $X$. 
\label{th:doublp_norm}
\end{lemma}
\begin{IEEEproof}
%
Let $r_i$ be the $i$'th row of $X$ and note $\| r_i \|_1 = 1$ for all $i \in \{1,\dots,n\}$, i.e., the sum of a row is equal to $1$.  Observe 
\[
\| X \|_F = \sqrt{\mathrm{Tr} \left( X X^*\right)} = \sqrt{\sum_{i=1}^n \| r_i\|_2^2 } \le \sqrt{\sum_{i=1}^n \| r_i\|_1^2 } \le \sqrt n,
\] 
where the first inequality follows because $\| \cdot \|_2 \le \| \cdot \|_1$.
\end{IEEEproof}

The second lemma establishes that $X^\star - X_{k}$ is a scaled doubly stochastic matrix.
\begin{lemma} \label{th:residual}
Let $X_{k} = \sum_{i=1}^k \theta_i P_i$ and suppose $X_k(a,b) \le X^\star(a,b)$ for all $a,b \in \{1,\dots,n\}$ and $k\ge 1$. Then,
\begin{itemize}
\item [\textup{(a)}] $\displaystyle \frac{X^\star - X_{k}}{1-{\sum_{i=1}^{k} \theta_{i}}}$ is doubly stochastic\\
\item [\textup{(b)}] 
$\displaystyle \| X^\star - X_k \|_F \le \sqrt{n} \left(1-\sum_{i=1}^{k} \theta_i \right)$
\end{itemize}
\end{lemma}
\begin{IEEEproof}
We start with (a). By assumption, $0 \le X_{k}(a,b)  \le X^\star(a,b) \le 1$ for all $a,b \in \{1,\dots,n\}$.  Hence, we only need to show that the sum of each row and column is equal to one. Observe
\begin{align*}
& \textstyle   (1- {\sum_{i=1}^{k} \theta_{i}} )^{-1} (X^\star - X_{k} ) \1 \\
& \quad  =  \textstyle (1-{\sum_{i=1}^{k} \theta_{i}})^{-1}  (X^\star - \sum_{i=1}^{k} \theta_i P_i ) \1 \\
& \quad =   \textstyle (1- {\sum_{i=1}^{k} \theta_{i}} )^{-1} ( X^\star \1 -  \sum_{i=1}^{k} \theta_i  P_i \1)  \\
& \quad =   \textstyle (1- {\sum_{i=1}^{k} \theta_{i}} )^{-1} (\1 -  \1 \sum_{i=1}^{k} \theta_i)  \\
& \quad =   \textstyle \1 (1- {\sum_{i=1}^{k} \theta_{i}} )^{-1} (1 - {\sum_{i=1}^{k} \theta_{i}})  \\
& \quad = \textstyle \1
\end{align*}
The same argument above can be used to show that $ (1- {\sum_{i=1}^{k} \theta_{i}} )^{-1} \1^T (X^\star - X_{k} ) = \1^T$, i.e., the sum of each column is equal to one. 

For (b), observe $\|  ( 1- \sum_{i=1}^{k} \theta_{i} )^{-1}  (X^\star - X_{k}) \|_F  = ( 1- \sum_{i=1}^{k} \theta_{i} )^{-1}  \| X^\star - X_{k} \|_F \le \sqrt n$ by Lemma \ref{th:doublp_norm}. Rearranging terms yields the result. 
\end{IEEEproof}

\subsection*{Proof of Lemma \ref{th:birkhoff_lemma}}

We start by proving the lower bound. We first note $\| X \|_F \ge 1$ for any doubly stochastic matrix $X$. Recall
\[
\| X \|_F \ge \| X \|_2 \coloneqq \sup \left\{ \frac{\| Xu \|_2 }{\| u \|_2} \ \text{with } u \in \R^n \ \text{s.t. } u \ne 0\right\}
\]
Let $u = \1$ in the equation above to obtain
\[
\| X \|_F \ge \frac{\|X \1 \|_2}{\| \1 \|_2} = \frac{\| \1 \|_2}{\| \1 \|_2} = 1,
\]
where $X \1 = \1$ follows since $X$ is doubly stochastic. Next, since $\frac{X^\star - X_{k}}{1-{\sum_{i=1}^{k} \theta_{i}}}$ is doubly stochastic by Lemma \ref{th:residual}, we have
\[
1 \le \left\| \frac{X^\star - X_{k}}{1-{\sum_{i=1}^{k} \theta_{i}}} \right\|_F  = \left( 1-{\sum_{i=1}^{k} \theta_{i}} \right)^{-1} \| {X^\star - X_{k}} \|_F
\] 
Rearranging terms yields the lower bound.

For the upper bound, observe 
\begin{align}
& \| X_{k} - X^\star \|_F^2 \notag \\
\text{(a)} &  = \| X_{k-1} + \theta_k P_k - X^\star \|_F^2  \notag\\ 
&  = \textstyle \| X_{k-1} - X^\star \|_F^2  + \theta^2_k \| P_k \|_F^2 \notag\\
& \quad \textstyle + 2 \theta_k \sum_{a,b} P_k(a,b) (X_{k-1}(a,b) - X^\star(a,b)) \notag \\
\text{(b)} &  \le \textstyle \| X_{k-1} - X^\star \|_F^2  + \theta^2_k \| P_k \|_F^2 - 2 \theta^2_k \sum_{a,b} P_k(a,b)^2  \notag \\ 
& = \textstyle \| X_{k-1} - X^\star \|_F^2  + \theta^2_k \| P_k \|_F^2 - 2 \theta^2_k n  \notag \\
\text{(c)} & \le \textstyle \| X_{k-1} - X^\star \|_F^2  + \theta^2_k n - 2 \theta^2_k n \notag \\
&  = \| X_{k-1} - X^\star \|_F^2  - \theta^2_k n \label{eq:thetaboundproofeq}
\end{align}
where (a) follows by Algorithm \ref{al:birkhoff_algorithm}, (b) by Eq.\ \eqref{eq:k1}, and (c) by Lemma \ref{th:doublp_norm}. 
Hence, 
\begin{align*}
\| X_{k} - X^\star \|_F^2 \le \left(1  - \frac{n \theta^2_k}{\| X_{k-1} - X^\star \|_F^2} \right) \| X_{k-1} - X^\star \|_F^2 
\end{align*}
Applying the argument recursively from $i=1,\dots, k$ 
\begin{align*}
\| X_{k} - X^\star \|_F^2  \le   \| X_0 - X^\star \|_F^2 \prod_{i=1}^k  \left(1  - \frac{n \theta^2_i}{\| X_{i-1} - X^\star \|_F^2} \right)
\end{align*}
Finally, since $X_0  =  \{ 0 \}^{n \times n} $ and $\| X^\star \|_F \le \sqrt n$ by Lemma \ref{th:doublp_norm}, 
\begin{align*}
\| X_{k} - X^\star \|_F^2  \le   n \prod_{i=1}^k  \left(1  - \frac{n\theta^2_i}{\| X_{i-1} - X^\star \|_F^2} \right)
\end{align*}
Taking square roots on both sides yields Eq.\ \eqref{eq:birkhoff_upper_bound}. 

To conclude, we show that $\theta_i \le  \frac{1}{\sqrt n} {\| X_{i-1} - X^\star \|_F}$ for all $i=1,2,\dots,k$. From Eq.\ \eqref{eq:thetaboundproofeq}, $0 \le \| X_{k-1} - X^\star \|_F^2  - \theta^2_k n $. Rearranging terms and taking square roots on both sides completes the proof.

\subsection*{Proof of Lemma \ref{th:lemma_decomposition}}

Since $ {n \theta^2_i}/{\| X_{i-1} - X^\star \|_F^2 } \ge \mu_{\min}$ by assumption, the upper bound in Lemma \ref{th:birkhoff_lemma} becomes $\| X_{k} - X^\star \|_F  \le  \sqrt n \left(1  -  \mu_{\min} \right)^{k/2}$.
Next, let $\epsilon = \| X_{k} - X^\star \|_F$ and write $\epsilon \le \sqrt n \left(1  -  \mu_{\min} \right)^{k/2} $. Rearranging terms yields 
\[
\left(\frac{1}{1- \mu_{\min}}\right)^{k/2} \le \frac{\sqrt n }{\epsilon}.
\]
Taking logs on both sides and further rearranging terms yields the result.

\subsection*{Proof of Theorem \ref{th:birkhoff_linear_convergence}}
We start by showing that we can design a subroutine \texttt{PERM} that returns a permutation with an associated weight that is uniformly lower bounded and satisfies the conditions in Eqs.\ \eqref{eq:k1}--\eqref{eq:k3}. We have the following lemma. 

\begin{lemma} Set $\mathcal I_k(\alpha)$ with $\alpha = \frac{1-\sum_{i=1}^{k-1} \theta_i}{(n-1)^{2}+1}$ is non-empty. 
\label{th:perm_set}
\end{lemma}
\begin{IEEEproof}
By Lemma \ref{th:residual},  $\frac{ X^\star - X_{k-1}  }{1-\sum_{i=1}^{k-1} \theta_i}$ is doubly stochastic, and so, by Carath\'eodory's theorem, we can write it as the convex combination of $(n-1)^2 + 1$ permutation matrices,  i.e., 
\[
 \frac{X^\star - X_{k-1} }{1-\sum_{i=1}^{k-1} \theta_i}  =  \sum_{j=1}^{(n-1)^2+1}  \beta_j P_j 
\] 
where $\beta_j \ge 0$ and $\sum_{j=1}^{(n-1)^2 + 1} \beta_j = 1$.
Next, note that since the permutation matrices and weights are non-negative, we have that
\begin{align}
\beta_j P_j (a,b) \le  \frac{X^\star(a,b) - X_{k-1}(a,b) }{1-\sum_{i=1}^{k-1} \theta_i} \label{eq:alpha_a}
\end{align}
holds for all $ a,b \in \{1,\dots,n\}$ and $j \in \{1,\dots,(n-1)^2+1\}$. Furthermore, since $\sum_{j=1}^{(n-1)^2 + 1} \beta_j = 1$, we have that
\begin{align}
\beta_j \ge \frac{1}{(n-1)^2 + 1} \label{eq:alpha_b}
\end{align}
for at least one $j \in \{1,\dots,(n-1)^2+1\}$. Let $\alpha = \beta_j$ such that the last equation holds. 
Combining Eq.\ \eqref{eq:alpha_a} and Eq.\ \eqref{eq:alpha_b}, we obtain that 
\begin{align*}
\frac{1-\sum_{i=1}^{k-1} \theta_i}{(n-1)^2 + 1} P (a,b) = \alpha P (a,b) \le X^\star(a,b) - X_{k-1}(a,b) 
\end{align*}
That is, there exists at least a permutation $P$ such that $X_{k-1} + \alpha P (a,b) \le X^\star(a,b)$, and so set $\mathcal I_k (\alpha)$ is non-empty.
\end{IEEEproof}

We are now in position to present the proof of Theorem \ref{th:birkhoff_linear_convergence}. By Lemma \ref{th:perm_set}, set $\mathcal I_k(\alpha')$  with $\alpha' = \frac{1-\sum_{i=1}^{k-1} \theta_i}{(n-1)^{2}+1}$ is non-empty. Now, observe that since $\alpha = \sqrt{\frac{\mu_{\min}}{ n} } (1-\sum_{i=1}^k \theta_i) \le \alpha'$ because $\mu_{\text{min}} = 1/n^3$, we have that $\mathcal I_k(\alpha') \subseteq \mathcal I_k(\alpha)$ and so $\mathcal I_k(\alpha)$ is non-empty. The rest of the proof follows as in Lemma \ref{th:lemma_decomposition} with $\mu_{\text{min}}  = \min_{i \in \{1,\dots,k\}} {n \theta^2_i}/{\| X_{i-1} - X^\star \|_F^2} $.

\subsection{Proofs of Section \ref{sec:frank-wolfe}}

\subsubsection*{Birkhoff polytope representation in vector form}

The Birkhoff polytope is the set that contains all doubly stochastic matrices. Recall we say that a nonnegative matrix is doubly stochastic if the sum of its rows and columns is equal to one. This corresponds to having $2n$ equality constraints. We can express these in vector form by defining matrices 
\begin{align*}
A' (in+1,in+j) &  = \begin{cases}
1 \qquad  i=0,\dots,n-1, \  j = 1,\dots,n\\
0 \qquad \text{otherwise}
\end{cases} \\
A''{(j+in,j)} &  = 
\begin{cases}
1  \qquad i=0,\dots,n-1, \  j = 1,\dots,n\\
0 \qquad \text{otherwise}
\end{cases}
\end{align*}
and then collecting them in $A= [A' ; A'']$. Next, define $ b \in \{ 1 \}^{2n}$. Any vector from $\R^d_+$ such that $Ax = b$ correspond to having doubly stochastic matrix in vector form.

For example, with $n=3$ we have
\begin{align*}A = 
\begin{bmatrix}
 1 & 1 & 1 & 0 & 0 & 0 & 0 & 0 & 0 \\
 0 & 0 & 0 & 1 & 1 & 1 & 0 & 0 & 0\\
 0 & 0 & 0 & 0 & 0 & 0 & 1 & 1 & 1\\
 1 & 0 & 0 & 1 & 0 & 0 & 1 & 0 & 0\\
 0 & 1 & 0 & 0 & 1 & 0 & 0 & 1 & 0\\
 0 & 0 & 1 & 0 & 0 & 1 & 0 & 0 & 1
\end{bmatrix}, \qquad
b = \begin{bmatrix}
1 \\ 
1\\
1\\
1\\
1\\
1
\end{bmatrix}.
\end{align*}

\subsection*{Proof of Lemma \ref{th:eccentricity}}

The eccentricity consists of two parameters. The diameter of the polytope ($M$) and its pyramidal width ($\lambda$). The diameter of the Birkhoff polytope is the maximum distance between two points in $\mathcal B$, which is the maximum distance between two vertices. Specifically, this is equal to $\|p -p' \|_2 = \sqrt{2n}$ where $p, p' \in \mathcal P$ are two vertices such that $p^T p' = 0$, i.e. have ones in different components.  

It is possible to obtain a lower bound on the pyramidal width of the Birkhoff polytope by using the fact that its extreme points are a subset of the extreme points of the unit cube in $d$ dimensions. Formally, $\mathcal P \subset \{0,1\}^d$ and so $ \mathrm{conv} (\mathcal P) \coloneqq \mathcal B \subset \mathcal C \coloneqq \mathrm{conv}( \{0,1\}^d)$. The latter means that the unit cube is ``extreme-point-wise denser'' than the Birkhoff polytope and so it has smaller pyramidal width. From Lemma 4 in \cite{LJ15} we can obtain that the pyramidal width of the Birkhoff polytope is lower bounded by $1 / \sqrt{d} = 1 / n$.\footnote{Recall that $d = n^2$.} Hence, $(\lambda /M)^2 \ge  {1}/{(2n^3)}$ as claimed. 

\subsection*{Proof of Theorem \ref{th:FCFW}}

This theorem is an application of  Theorem 1 in \cite{LJ15} with the quadratic objective function $f(x) = (1/2)\| x - x^\star  \|_2^2$ and set $\mathcal B$. This theorem says that
\[
 \| x_k - x^\star \|_2^2 \le \| x_0 - x^\star \|_2^2 \exp \left( - \frac{\mu}{4L} \left(\frac{\lambda}{M}\right)^2 k \right)
\]
The term $\|  x_0 - x^\star\|_2$ can be upper bounded by $\sqrt{2n}$, which is the maximum Euclidean distance between two points in $\mathcal B$ (see the proof of Lemma \ref{th:eccentricity}). The condition number $\mu/L$ is equal to $1$ because the objective function is quadratic and $({\lambda}/{M})^2 \ge 1/(2n^3)$ by Lemma \ref{th:eccentricity}. Hence,
\[
 \|  x_k - x^\star \|_2^2 \le 2n  \exp \left( - \frac{k}{8n^3} \right) .
\]
To conclude, let $\epsilon^2  =  \|  x_k - x^\star\|_2^2$ and write $\epsilon^2 \le 2n  \exp ( - {k}/({8n^3}) )$. By expressing $k$ as a function of $\epsilon$ in the last equation, we obtain the stated result.

\end{document}